\documentclass[12pt,a4paper]{article}
 \usepackage{amsmath,amstext,amssymb,amscd}
 \usepackage{graphicx}
\begin{document}

\title{Complementary choice functions }
\author{Danilov V.I. \thanks{Central Institute of Economics and Mathematics of the RAS, Nahimovskii prospect, 47, 117418, Moscow; email: vdanilov43@mail.ru.}}
\date{\today}
\maketitle

\begin{abstract}
The paper studies complementary choice functions, i.e. monotonic and consistent choice functions. Such choice functions were introduced and used in the work \cite{RY} for investigation of matchings with complementary contracts. Three (universal) ways of constructing such functions are given: through pre-topologies, as direct images of completely complementary (or pre-ordered) choice functions, and with the help of supermodular set-functions.


\end{abstract}

Choice functions have long attracted the attention of economists. First, to describe the demand and generally rational behavior. Later, the choice functions were used in the theory of stable partnerships and contracts. A special role plays the so-called substitute choice functions, introduced (for another reason) by Plott \cite{P}. In 2020, an interesting article by Rostek and Yoder \cite{RY} has appeared, in which the attention was attracted to the opposite property of complementarity. They shown that this property also guarantees the existence of stable contract systems. Surprisingly, this property did not attract attention of workers on the choice theory. In this note, we want to to fill this gap and focus on the mathematical side of the matter - how to set such choice functions and how they are related to supermodularity property.

It should be said that supermodularity (and the dual concept submodularity) has long attracted the attention of specialists in combinatorial mathematics. In particular, it received a great response
a polynomial algorithm for finding of the maximum of such functions (see \cite{Schr}). They work with  the lattice $2^X$ of subsets of some basic set $X$ (usually finite). Topkis \cite{Top} developed a theory of supermodular functions for arbitrary lattices. As a rule, if not touch on subtle issues and sit at the level of definitions, work with arbitrary lattices does not cause difficulties. In order not to complicate the summary, we will work with the lattice of subsets, and only in the final section will talk about arbitrary lattices.

\section{Definitions and examples}

      Let $X$ be a set. A \emph{choice function} (CF) on $X$ is a mapping
                                                         $$
                                 f: 2^X \to 2^X,
                                                           $$
such that $f(A)\subseteq A$ for any $A\subseteq X$. There are very a lot of  choice functions, and we will be interested only in those that reflect the choice of "best"\ elements in a sense. One might think that CF formalizes the behavior of some decision-maker: having access to a set of items $A\subseteq X$ ("menu"), the person chooses subset $f(A)\subseteq A$ of its most interesting subjects. Such interpretation justifies attention to the following property \emph{consistency} of CFs (other names are also used: the outcast condition, the condition of independence from rejected
alternatives, IRA, etc.)
$$
      \textrm{if }f(A)\subseteq B\subseteq A,\textrm{ then }f(B)=f(A).
                                          $$
In other words, the choice does not change if we remove some rejected elements. Usually this common condition of "rationality"\ often augmented by additional conditions reflecting additional presentations about the "best". We will limit ourselves to two antipodal structural conditions. \medskip

\textbf{Definition.} A CF $f$ is called \emph{subadditive} (resp., \emph{superadditive}) if the following inclusions are fulfilled for any $A,B\subseteq X$
                                                       $$
f(A\cup B)\subseteq      f(A)\cup f(B)
      $$
(resp., the opposite inclusions $f(A)\cup f(B)\subseteq f(A\cup B)$).\medskip

           Subadditivity is better known in the following formulation:
      $$
                             \textrm{if }A\subseteq B,\textrm{ then  }     f(B)\cap A\subseteq f(A).
        $$
In other words, if some element from the set $A$ is chosen in a larger set $B$, then it is chosen in the smaller set $A$. For this reason Aizerman and Malishevski \cite{AM} called this property as the heredity property. How Blair \cite{Bl} rightly points out, this property expresses that the elements of $f(A)$ are chosen because of their own qualities, not because of their connection with other elements ("members of $C_i(A)$ are wanted for their own sake, not because of potential benefits from interaction with other members"). For this reason, they are called substitutable.  CFs which satisfy this property (especially together with the consistency), have intensively studied in the literature (\cite{P,AM, DK}) and found an application in the theory of stable contracts and in non-monotonic logic.

Superadditivity, as it is easy to understand, is equivalent to the simpler monotonicity property:
$$
                              \textrm{if }A\subseteq B,\textrm{ then }f(A)\subseteq f(B).
                                               $$

\textbf{Definition.} A CF $f$ is called \emph{complementary} if it satisfies the conditions of consistency and monotonicity.\medskip

Such CFs did not attract an attention until the appearance of the work of Rostek and Yoder \cite{RY} about complementary contracts. In the economic literature devoted to demand, theories often divide goods into substitute and complementary ones. Condition heredity (=subad\-ditivity) reflects substitutability, whereas the monotonicity condition (=superad\-ditivity) reflects complementarity. In the theory of stable contracts, it has long been understood, that the condition of substitutability guarantees the existence of stable contract systems. From the filing of Milgrom, the belief prevailed that the substitutability condition is even necessary for existence. Maybe this       explains the one-sided attention to substitutability and neglect of the complementarity. We tried to challenge this misconception (\cite{DKL}), but in vain. The work of Rostek and Yoder opened the way to study of complementary CFs.\medskip

\textbf{Lemma 1.} \emph{Let $f$ be a monotone CF. It is complementary if and only is it is       idempotent, that is, if $f(f(A))=f(A)$ for any $A$.}\medskip

Proof. Let $f$ be consistent. Since $f(A)\subseteq f(A)\subseteq A$, then we conclude that $f(f(A))=f(A)$.

Conversely, let $f$ be idempotent. And let $f(A)\subseteq B\subseteq A$. In force of the       monotonicity, from  $f(A)\subseteq B$ we have $f(A)=f(f(A))\subseteq f(B)$. On the other hand, from $B\subseteq A$ we get $f(B)\subseteq f(A)$, from where the equality $f(A)=f(B)$ holds. $\Box$\medskip

Before proceeding to the study of complementary CFs, it is useful to give some examples of such CFs. Especially since they are not only particular examples, but are building blocks for arbitrary complementary CFs.\medskip

\textbf{Example 1.} Let $K$ be a subset of $X$. Define the CF $f_K$ as follows: $f_K(A)=K$ if $K\subseteq A$, and $f_K(A)=\emptyset$ in other cases. It is easy to understand that the CF $f_K$ is complementary. This CF reflects the DM's interest exclusively in the set $K$: if set $K$ is available, $K$ is selected; otherwise nothing is selected. For this reason, these CFs can be called \emph{packaged}.\medskip

\textbf{Example 2.} Let $\le$ be some preorder on $X$. An (ordinal) ideal is a subset $I$ in $X$, which with each element contains smaller ones, $x\le y$ and $y\in I$ implies $x\in I$. Define the CF $f_\le$ as follows: $f_\le (A)$ is equal to the largest ideal contained in $A$. Again       obviously, the CF $f_\le$ is complementary.\medskip

\textbf{Example 3.} Let $k$ be a natural number or $\infty$. Define $f_k$ so: $f_k(A)=A$ if $A$ contains $k$ or more elements, and $f(A)=\emptyset $ if $|A|<k$. This CF is also complementary.\medskip

\textbf{Example 4.} A variant of the previous one. Let $f(A)=A$ if the set $X-A$ is finite, and $f(A)=\emptyset$ otherwise.

           \section{The structure of complementary CFs}

Let CCF($X$) denote the set of complementary CFs on the set $X$. This set has a natural operation of union. If $(f_i, i\in I)$ is a family of CFs, then CF $\cup _if_i$ is given by the obvious way:
$$
                              (\cup _if_i)(A)=\cup _if_i(A).
                                                                  $$

\textbf{Proposition 1.} \emph{If all $f_i$ are complementary CFs, then $f=\cup _if_i$ is a       complementary CF.}\medskip

Proof. Since $f$ is obviously monotonic, then by virtue of Lemma 1 it is enough to check the idempotency of $f$. Due to monotonicity of $f$, the set $f(f(A))=f(\cup _if_i(A))$ contains $f(f_j(A))=\cup      _if_i(f_j(A))$ for any $j$ and, in particular, $f_j(f_j(A))=f_j(A)$. Therefore $f(f(A))$ contains all $f_j(A)$, that is, it contains $\cup _if_i(A)=f(A)$. And since $f(f(A)) \subseteq f(A)$, we get the equality $f(f(A))=f(A)$. $\Box$\medskip

Thus, the poset CCF($X$) is an upper semilattice regarding the $\cup$ operation. This hints at the fact that any complementary CF can be represented as a union of elementary (join-irreducible) complementary CFs. Note that the packaged CFs from Example 1 are just join-irreducible. We shall show that any complementary CF is indeed represented as a join (a union) of packaged CFs. To do this, we will need the following\medskip

\textbf{Definition.} A subset of $A\subseteq X$ is called \emph{open} (for CF $f$) if $f(A)=A$.\medskip

For example, the empty set $\emptyset$ is open. The set of open subsets in $X$ is denoted as $St(f)$. In the case of complementary CF this set coincides with the image $f$, $St(f)=Im(f)$.\medskip

\textbf{Theorem 1.} \emph{Let $f$ be a complementary CF. Then $f$ is represented as the union of elementary packaged CFs $f_K$, where $K$ runs through the set       $St(f)$,}
      $$
                                f=\cup _{K\in St(f)} f_K.
                                $$

Simply put, $f(A)$ is the largest open subset of $A$, that is, the `interiority' of $A$.\medskip

Proof. Obviously, $f_K\subseteq f$ for any open $K$. So that it remains to check the reverse inclusion of $f\subseteq\cup _{K\in St(f)}
f_K$. Let $A$ be
an arbitrary subset in $X$, and $K=f(A)$. Then $f(A)=f_K(A)$. $\Box$\medskip

           This forces us to pay an attention to the set of $St(f)$.\medskip

\textbf{Lemma 2.} \emph{Let  $f$ be a monotonic CF. If $(S_i, i\in I)$ is a family of open subsets, then $S=\cup _i S_i$ is also open.}\medskip

In fact, due to monotonicity, $f(S)$ contains $f(S_i)=S_i$ for any $i\in I$, so that $f(S)$ contains $S$, from where $f(S)=S$. $\Box$\medskip

Thus, the set $St(f)$ is closed with respect to any unions. In \cite{DK} such systems of sets were called pre-topologies. The only difference from topology is that we do not require the openness of the intersection of open subsets. The set $f(A)$ should be understood as the `interiority' of $A$.

Conversely, let $\mathcal P\subseteq 2^X$ be some pre-topology on $X$. Define CF $f_\mathcal P $ as the union of the packaged CFs $f_K$, where $K$ runs through $\mathcal P$, that is open sets of the pre-topology. Then $f_\mathcal P$ is a complementary CF, moreover, as it is easy to check, $St(f_\mathcal P )=\mathcal P $. Thus we get that the set (and even the poset) CCF($X$) of complementary CFs on $X$ can be identified with the set of pre-topologies on $X$.

The use of pre-topologies is useful in the following respect. The fact is that there is a simple way to set a pre-topology: to do this, you can take any system (`a base') $\mathcal B$ of subsets in $X$ and add to it all possible unions of sets from $\mathcal B$. In terms of the base $\mathcal B $ the corresponding complementary CF $f$ is arranged as follows: $f(A)$ is the union of those $B\in\mathcal B$ which contained in $A$,
$$
      f(A)=\cup _{B\in \mathcal B , B\subseteq A} B.
$$

      \section{Continuous and completely complementary CFs}

      Let's continue to use the topological terminology.

Let $f$ be a complementary CF on $X$, and $x\in X$. Let's say that $N\subseteq X$ is a \emph{neighborhood} of the point $x$ if $x\in f(N)$. Such terminology is justified by that $x$ belongs to any of its neighborhoods, and any superset of a neighborhood is also a neighborhood. Denote by $\mathcal N(x)$ the set of neighborhoods of the point $x$. In these terms, CF $f$ can be set by the compact (tautological) formula
                                    $$
                                f(A)=\{a\in A,\ A\in      \mathcal N (x)\}.
                                      $$
Instead of arbitrary neighborhoods, one can use open neighborhoods, that is such neighborhoods of a point $x$ which are open. Let's denote the set of open neighborhoods as $\mathcal U(x)$. Each neighborhood $N$ of a point $x$ contains an open neighborhood, namely $f(N)$. Therefore, the previous the formula can be rewritten as
                                                    $$
      f(A)=\{a\in A,\textrm{ there exists } U\in \mathcal U(a)\textrm{ such that } U\subseteq A\}.
                   $$

From this point of view, minimal neighborhoods (which are obviously open) are of interest. In general, they may not to exist. Such an example gives CF $f_\infty$ from Example 3 with infinite $X$.\medskip

\textbf{Definition.} A CF $f$ is \emph{continuous} at a point $x\in X$ if any neighborhood of $x$ contains a minimal neighborhood of $x$. A CF is \emph{continuous} if it is continuous at all points.\medskip

For example, if $X$ is a finite set, then any CF is continuous. Denote by $\mathcal M(x)$ the set of minimal neighborhoods of a point $x\in X$. Note that this family of subsets satisfies the following three properties a)-c):

a) $x\in S$ for any $S\in\mathcal M(x)$,

b) subsets of $\mathcal M(x)$ are incomparable by inclusion,

c) if $y\in S\in\mathcal M(x)$, then there exists $T\in\mathcal M(y)$, such that $T\subseteq S$.

It is clear that for continuous CF $f$ the previous formula can be corrected:
                                                              $$
      f(A)=\{a\in A,\textrm{ there exists }S\in\mathcal M(a),\ S\subseteq       A\}.
                                                                $$
This suggests a way to build such CFs. Suppose that, for any points $x\in X$, a set of $\mathcal M(x)$ subsets is given, so that the properties a)-c) are satisfied. Define a CF $f$ by the above formula
$$
      f(A)=\{a\in A,\textrm{ such that there exists }S\in \mathcal M(a),\ S\subseteq A\}.
                   $$
Such a CF $f$ is obviously monotonous. By virtue of c) it is idempotent, and by Lemma 1
is complementary. A subset $N$ is a neighborhood of $x$ if and only if $N$ contains some $S$ of $\mathcal M(x)$. Therefore (by virtue of b)) set $\mathcal M(x)$ is exactly the set of minimal neighborhoods of a point $x$ for the CF $f$.

Especially interesting is the case when, for any point $x$ of $X$, the set $\mathcal M(x)$
consists of a single set $S(x)$ (containing, recall, $x$). In this case
$$
      f(A)=\{a\in A, \ S(a)\subseteq       A\}.
        $$
Of course, $S(x)$ is open for $f$. At the same time, the intersection any open subsets are open. Indeed,  let $(S_i, i\in I)$ be an arbitrary family of open sets, and $S=\cap _i S_i$. If  $x\in S$, then $x$ will give any $S_i=f(S_i)$, which means that $S(x)\subseteq S_i$ for any $i\in I$ and $S(x)\subseteq S$, so that $x\in f(S)$. And we get that $S\subseteq f(S)$, that is $S= f(S)$. \medskip

\textbf{Definition.} Say that a consistent CF $f$ is \emph{completely complementary} if, for any family $(A_i, i\in I)$ menu $A_i$,
$$
                              f(\cap _i A_i)=\cap _i f(A_i).
                                         $$
      Note that such a CF is monotonic and, therefore, is complementary.\medskip

An example of completely complementary CF gives Example 2. For CF from Example 2, the set      $S(x)$ is the ideal generated by $x$,  $S(x)=\{y,\ y\le x\}$.\medskip

\textbf{Proposition 2.} \emph{For a consistent CF $f$ the following assertions are equivalent:}

1) $f$ \emph{is given as in Example 2};

2) $f$ \emph{is completely complementary};

3) $f$ \emph{is continuous and $\mathcal M(x)$ consists of a single set.}\medskip

Proof. Obviously, that CFs from Example 2 are completely complementary. If $f$ is completely complete, then the intersection of all the neighborhoods of the point $x$ is the minimal neighborhood of $x$, so that 2) implies 3). Finally, let 3) be fulfilled, and $S(x)$ is the only minimal neighborhood of $x$. Let us define the binary relation $\le$ on $X$, assuming $y\le x$ if and only if $y\in S(x)$. This relation is obviously reflexive. Transitivity follows from
the condition c). $\Box$\medskip

Completely  complementary CFs are also universal in a sense. That is, they allow to build any complementary CF using the operation of direct image generalizing (in some meaning) the union operation. Namely, imagine that we are given a mapping of sets $\varphi :Y\to X$, and CF $g$ on $Y$. Then we can define (see \cite{DK}) CF $f=\varphi_*(g)$ by $X$, assuming (for $A\subseteq X$)
        $$
      f(A)=\varphi       (g(\varphi ^{-1}(A))).
          $$

\textbf{Proposition 3.} \emph{If the CF $g$ is complementary, then $f=\varphi_*(g)$ is complementary as well.}\medskip

Proof. Obviously, $f$ is monotonic. Its consistency is proven in \cite{DK}, Proposition 8. One can, however, not refer to \cite{DK}, but directly check the idempotency of $f$. $f(f(A))=\varphi (g(\varphi ^{-1}(\varphi (g(\varphi ^{-1}(A)))))$. Put $B=g(\varphi^{-1}(A))$, so that $ff(A)=\varphi(g(\varphi ^{-1}\varphi(B)))$, whereas $f(A)=\varphi (B)$. Since $B\subseteq \varphi ^{-1}\varphi(B)$, then (by virtue of monotonicity and idempotence of $g$) $B=g(B)\subseteq g(\varphi ^{-1}\varphi (B))$. Therefore
                                                                        $$
      f(A)=\varphi (B)\subseteq       \varphi g\varphi ^{-1}\varphi (B)=ff(A).
                                                             $$

Thus, the operation of direct image allows to build complementary CFs. At the same time, we can assume that the CF $g$ on $Y$ is completely complementary.\medskip

\textbf{Theorem 2.} \emph{Let $f$ be a complementary CF on $X$. There are: a set $Y$ together with a mapping $\varphi :Y\to X$ and a completely complementary CF $g$ on $Y$, such that $f=\varphi _*(g)$.}\medskip

Proof. We will explicitly construct such $Y$, $\varphi$ and $g$. Let $\mathcal U$ be a system of open neighborhoods for CF $f$, so that (for any menu $A$)
$$
f(A)=\{a\in A,\textrm{ there exists }U\in \mathcal U(a),\textrm{ such that }U\subseteq A\}.
                                             $$
Form the set $Y$ as the set of pairs $(x,U)$, where $U$ is open set, and $x\in U$. So in fact $Y=\coprod_x\mathcal U(x)$. The mapping of $\varphi$ will be set naturally: $\varphi(x,U)=x$. Finally, consider the following binary relation $\le$ on  $Y$, assuming $(x,U)\ge(y,V)$ if $V\subseteq U$. This relation is reflexive and transitive. Define the CF $g$ on $Y$ as the CF associated with this preorder $\le$ (as in Example 2). That is, for any $B\subseteq Y$
                                                                $$
      g(B)=\{b\in B, \textrm{ such that everything less than } b \textrm{ is contained in } B\}
                                                                 $$
(that is $b'\le b$ implies $b'\in B$). In other words, if $b=(x,U)$, then for any open $V$ lying in $U$, any pair of $(y,V)$ ($y\in V$) must also belong to $B$. This CF is obviously (by virtue of Proposition 2) completely complementary. And it remains to check that $\varphi _*(g)=f$. That is, for $A\subseteq X$  it is true that
      $$
      f(A)=\varphi       g(\varphi ^{-1}(A)).
        $$

Inclusion $\subseteq$. Let $a\in f(A)$. Consider the pair $b=(a,U)$, where $U=f(A)$. Obviously, this pair is selected (FB $g$) from the set $B=\varphi ^{-1}(A)$. In fact, let $b'=(x,V)\le b$, that is, $x\in V\subseteq U=f(A)\subseteq A$. Since $\varphi(b')=x$ is contained in $A$, then $b'\in
      \varphi ^{-1}(A)$.

Inclusion $\supseteq$. Let $a\in\varphi g(\varphi^{-1}(A))$, then there is for some an open $U$ pair $(a,U)$ is selected from the set $\varphi^{-1}(A)$. It means that, for any smaller pair $(x,V)$ ($x\in V\subseteq U$), we have $(x,V)\in \varphi ^{-1}(A)$, that is $x\in A$. In particular, for any $x$ of $U$ we have $x\in A$, that is $U\subseteq A$. But then $a\in U=f(U)\subseteq f(A)$. $\Box$\medskip

\textbf{Remark 1.} In the case when CF $f$ is continuous, the construction can be made more economically. Namely, build $Y$ as a set of pairs $(x,U)$, where $U$ is a minimal neighborhood of the point $x$.\medskip

\textbf{Remark 2.} Although complementary CFs are somewhat antipodal to substitutable ones, there is clearly a certain parallel. This is the presentation via join-irreducible CFs (Theorem 1), and the "simplification"\ by the operation of direct image (Theorem 2).

      \section{Supermodular functions}

There is another interesting way to generate complementary CFs, using the so-called supermodular set-functions. It was in this way Rostek and Yoder showed the complementarity of CFs they were interested in. Here we shall assume that $X$\emph{ is a finite set}.

Recall that if $u:2^X\to\mathbb R$ is a set-function on $X$, then it generates CF $f=f_u$ by the following rule: $f(A)$ is the subset of $A$ with the largest value of $u$. Here we need to assume,      what such a subset is unique (or somehow stands out among the largest ones); more will be said below.\medskip

\textbf{Definition.} A set-function $u:2^X\to\mathbb R$ is called \emph{supermodular} if (for any $A$ and $B$ in $X$) the following inequality are satisfied
                     $$
       u(A)+u(B)\le u(A\cap B)+u(A\cup B).
                       $$
If the inequalities are opposite, $u$ is called \emph{submodular}. If both, $u$ is called \emph{modular}. As a rule, it is assumed, that $u(\emptyset )=0$.\medskip

Such set-functions have been intensively studied and used in combinatorics. In discrete mathematics (in the graph theory, see, for example, \cite{LP}), in the theory of measure (Choquet capacity), in the theory of cooperative games (Shapley's convex games), in the theory of stable assignments  (\cite{Sher}), in economics when studying comparative statics (\cite{MiSh, Top}). As Rostek and Yoder noticed, such set-functions are related to monotone CF (a fact not previously noticed only       due to lack of interest in monotone CFs).\medskip

           \textbf{Lemma 3.} \emph{Let $u$ be a supermodular set-function on $X$. The set its
      maxima (more precisely, the set $Argmax(u)$) forms a sublattice in $2^X$.}\medskip

           Obviously.\medskip

With such a function $u$, we will associate CF $f=f_u$, putting $f(A)$ to be the least (by inclusion) subset of $A$ with the largest value of $u$.\medskip

\textbf{Proposition  4} (Rostek and Yoder). \emph{For a supermodular set-function $u$ the CF
      $f=f_u$ is monotonic (and thus is complementary).}\medskip

Proof. Let $A\subseteq B$; we need to show that $A^*=f(A)\subseteq f(B)=B^*$. Apply the supermodularity inequality to $A^*$ and $B^*$:
                                                            $$
          u(A^*)+u(B^*)\le u(A^*\cap B^*)+u(A^*\cup B^*).
                                                              $$
Since $A^*\cup B^*$ is contained in $B$, then $u(A^*\cup B^*)\le u(B^*)$. Similarly, $u(A^*\cap B^*)\le u(A^*)$. So, there is equality here and there. In particular, $u(A^*\cap B^*)=u(A^*)$. And from the minimality of $A^*$ we get $A^*\subseteq A^*\cap B^*$, that is, $A^*\subseteq B^*$. $\Box$ \medskip

Thus, any supermodular set-function gives a complementary CF. This raises two questions. The first one - where to get supermodular functions from? how to build them? And the second one - do them give ALL the complementary CFs? We will show below that the answer to the second question is affirmative. But for this we will still have to build supermodular set-functions. So let's say briefly about the well-known construction of some supermodular function.

A simple but very useful remark is that the sum of supermodular functions is again supermodular. So what remains to look for "basic"\ (or extreme) supermodular set functions. It seems that the complete list of such functions is unknown, but some subset can be offered. Namely, let $U$ be a subset of $X$. We form the following "elementary"\ set-function $e_U$: for $A\subseteq X$
                          $$
      e_U(A)=1\textrm{ if } U\subseteq A, \textrm{ and $=0$ otherwise}.
                            $$

      \textbf{Lemma 4.} \emph{This function is supermodular.}\medskip

In fact, let $A$ and $B$ be given; we need to check the supermodularity inequality. If $e(A)=e(B)=0$, then the inequality is obvious. Let $e(A)=0$, and $e(B)=1$. The second means that $U\subseteq B$. But then $U\subseteq A\cup B$, so $e(A\cup B)=1$. Finally, let $e(A)=e(B)=1$. Then $A$ and $B$ contain $U$, from where $A\cap B$ and $A\cup B$ contain $U$ as well and the inequality is satisfied again (as equality). $\Box$\medskip

\textbf{Theorem 3.} \emph{Let $f$ be a complementary CF on $X$. Suppose that the set X is finite. Then there is a supermodular set-function $u$, which generates $f$.}\medskip

Proof. Here it will be more convenient to work not with CF $f$, but with the corresponding pre-topology $\mathcal U$ (consisting, as mentioned earlier, of open subsets). Let $U\in\mathcal U$ (that is f(U)=U). Above, we have already introduced the "elementary"\ set-function $e_U$ and
checked its supermodularity.

Now let us put $u=\sum _{U\in \mathcal U} e_U$. Actually, $u(A)$ is equal to the number of open subsets in $A$. As the sum of supermodular functions, it is also supermodular. We claim that      the Cf $f_u$, generated $u$, coincides with $f$. In fact, let $A$ be arbitrary menu, and $A^*=f(A)$ be the largest open subset in $A$. It is clear that the maximum value of the function is $u$ (on subsets $B$ contained in $A$) takes exactly on those $B$ that contain $A^*$. And the smallest among them is just $A^*$. $\Box$\medskip

Do not think that this simple reasoning is perfect original. In a similar style in graph theory one establishes that the number of edges leaving a subset of vertices $A$ depends submodularly on $A$ (see also the example of hypergraph cut functions from 44.1a in \cite{Schr}).\medskip

\textbf{Remark.} The set-function from theorem 3 is even a monotone supermodular function. Such a function together with each maximal $A$ contains and all larger sets. The question arise --
is there a supermodular function $u$ that (given a set $A$) takes on $2^A$ the only one maximum (equal to $f(A)$)? The answer is yes. To do this, take instead of the function, proposed in the proof, the function $u_\varepsilon$, given by the equality
                                                             $$
                                 u_\varepsilon (B)=u(B)-\varepsilon |B|.
                                                               $$
It is also supermodular (due to the modularity of the function $B\mapsto |B|$). And with
a small $\varepsilon >0$, it reaches the maximum at $A$.\medskip

Theorem 3 and the antipodality of substitution and complementarity suggest that there should be a "dual"\ statement: a CF, generated by a submodular set-function $u$, is substitute. Surprisingly, this is not the case.\medskip

\textbf{Example 5.} Let $X=\{a,b,c\}$, and the set function u is given by the following table:

\begin{center}
\begin{tabular}{| c|c |c | c|c |c |c |c |c |}
  \hline
  % after \\: \hline or \cline{col1-col2} \cline{col3-col4} ...
  $A$ & $\emptyset$ & $a$ & $b$ & $c$ & $ab$ & $ac$ & $bc$ & $abc$ \\ \hline
  $u(A)$ & 0 & 3 & 2 & 2 & 2 & 2 & 4 & 1 \\   \hline
\end{tabular}
\end{center}

      It is directly verified that $u$ is submodular. For example,
                                          $$
                         u(ab)+u(ac)=2+2\ge u(a)+u(abc)=3+1.
                                            $$
Let CF $f$ be generated by the function $u$. It is clear that $f(X)=\{b,c\}$, whereas     $f(\{a,b\})=\{a\}$. We see that $b$ is chosen in a larger $X$, but is not chosen in      lesser $\{a,b\}$. So $f$ is not substutute.

This is quite strange, because for modular-set functions everything is OK: the corresponding CF $f$ is "conditionally constant"\ (that is, it is substitute and complementary).\medskip

Different supermodular functions $u$ can generate one and the same complementary CF $f$. Recall that $f(A)$ is defined as the (minimal "by inclusion") subset in $A$ with the largest value of $u$. This means that it is not so much the function $u$ matters as the order relation generated by it on $2^X$. Therefore, it is natural to turn to formulations in terms of pre-orders on $2^X$.\medskip

\textbf{Definition.} A weak order (i.e. complete preorder) $\preceq$ on the set $2^X$ is called  \emph{supermodular}\footnote{In the literature usually the term "quasi-supermodularity"\! is used applying to the function, setting this order.} if (for any $A$ and $B$ in $X$) either
$A\preceq A\cap B$, or $B\preceq A\cup B$, and if $A\cap B \prec A$, then $B\prec A\cup B$.\medskip

Let $A$ and $B$ be maximal (hence equivalent) subsets in $X$ with respect to a supermodular weak order $\preceq$. Then $A\cap B $ also is maximal. In fact, let it be less than $A$, $A\cap B \prec A$. Then $B\prec A\cup B$, contrary to the maximality of $B$. In other words, there is the only minimal (by the inclusion) maximal (relative to order $\preceq$) subset. The same is true for any subset in $X$. Thus, for any $A\subseteq X$ there exists $A^*$ such that

           a)  $A^*\subseteq A$,

           b) $B\subseteq A$ implies $B\preceq A^*$,

c) if $A^*\preceq B \subseteq A$, then $A^*\subseteq B$.

Obviously, such $A^*$ is unique. In fact, let $A^{**}$ satisfies the same a)-c). Assuming in c) $B=A^{**}$, we get that $A^*\subseteq A^{**}$. But symmetrically $A^{**}\subseteq A^*$ whence equality.

Let us define the CF $f$ on $X$ assuming $f(A)=A^*$.\medskip

      \textbf{Proposition 6.} \emph{The CF $f$ is complementary.}\medskip

Proof. Obviously, it is consistent. But it also is monotone. Reasoning as above. Let $A\subseteq B$; we need to show that $A^*\le B^*$.

Consider $A^*\cap B^*$; since this set is in $A^*$, it is $\preceq A^*$. If it is equivalent to $A^*$, then by virtue of the minimality of $A^*$ we get the equality $A^*\cap B^*=A^*$, that is, $A^*\le B^*$. If $A^*\cap B^*\prec A^*$, then  $B^*\prec A^*\cup B^*$, which contradicts the maximality of $B^*$ in $B$. $\Box$

      \section{Generalization to lattices}

The previous concepts and reasoning were so primitive and hinted that they should work in a more general situation.  See also \cite{CMY} who gave examples of interesting lattices for economical application. Moreover, the concept of supermodularity in the book \cite{Top} was introduced for       arbitrary lattices. Monotonicity  and complementarity also do not meet difficulties.

Let $L=(L,\le )$ be a complete (for example, finite) lattice. Three examples may be useful: distributive lattices, products of chains, and finally the Boolean lattice $2^X$.

A choice function on $L$ is a contracting mapping $f:L\to L$, that is, $f(x)\le x$ for any $x\in L$. A CF $f$ is monotonic if $x\le y$ implies $f(x)\le f(y)$. A CF $f$ is consistent if $f(x)\le y\le x$ implies $f(y)=f(x)$. As before, we call A CF $f$ complementary if it is monotonous and consistent. Almost all facts and reasonings established earlier for the lattice $2^X$ are carried over for the case of arbitrary lattices. We will not repeat the proofs.

To understand the structure of complementary CFs, the concept of fixed elements is again important. An element $x\in L$ is fixed (relative to CF $f$) if $f(x)=x$. Applying the consistency condition to the chain $f(x)\le f(x)\le x$, we get the equality $f(f(x))=f(x)$. So $f(x)$ is fixed for any $x\in L$. The set $Fix(f)=f(L)$ of fixed elements is (as before) a complete join-semilattice in $L$, that is, for any family $(x_i, i\in I)$ of fixed elements, their join $\vee_i x_i$ is a fixed element.\medskip

\textbf{Proposition 7.} 1) \emph{A complementary CF $f$ is uniquely determined by the set $Fix(f)$.}

2) \emph{If $F$ is a complete join-semilattice in $L$, then there exists a (unique) complementary CF $f$ such that $F =Fix(f)$.}\medskip

Of course, $f(x)=\bigvee_{z\in F, z\le x} z$.\medskip

Now let's go towards the supermodularity (see \cite{Top}). Supermodularity of a function $u$ on the lattice $L$ is actually defined as above:
                                $$
      u(x)+u(y)\le u(x\wedge y)+u(x\vee       y)
                                       $$
for any $x,y\in L$ (here $\le$ is the order relation in $\mathbb R$). As before, the set of elements in $L$ in which u reaches the maximum value (i.e. $Argmax(u)$), is a sublattice in $L$. In particular, if $L$ is a finite lattice then there is a minimal (in the sense of the order $\le$ of the lattice $L$) element in $Argmax(u)$, which we denote as $f(1)$.

In fact, all of this applies to the $u$ constraint on the sublattice $L(x)=\{y\in L, y\le x\}$; the minimal element in $Argmax(u|L(x))$ is denoted as $f(x)$. Thus we get CF $f=f_u$ by $L$. As before this CF is complementary.

Thus, every supermodular function on the lattice $L$ gives a complementary CF.  A more fundamental fact is that any complementary CF can be obtained in this way.\medskip

\textbf{Theorem 4.} \emph{Let $f$ be a complementary CF on a finite lattice $L$. Then there is a supermodular function $u$ on $L$ (and even monotone and with integer value), which generates $f$.}\medskip

The proof is the same as the proof of Theorem 3. Define $u(x)$ to be equal to the cardinality of the set $Fix(f)\cap L(x)$.

      \end{document}